\documentclass[12 pt]{article}

\usepackage{CJK}
\usepackage{amssymb}
\usepackage{bbm}
\usepackage{amsfonts}
\usepackage{mathrsfs}
\usepackage{graphicx}
\usepackage{amsmath}
\usepackage{amsthm}
\usepackage{xypic}
\usepackage{graphicx}
\usepackage{times}
\usepackage{geometry}
\usepackage{color}
\usepackage{indentfirst}
\usepackage{cases}
\usepackage{hyperref}
\usepackage{graphicx}
\usepackage{float}
\usepackage{amsthm}
\usepackage{amssymb}
\usepackage{amsmath}
\usepackage{mathrsfs}
\usepackage{indentfirst}
\usepackage{geometry}
\usepackage{authblk}
\usepackage{hyperref}

\hypersetup{
    colorlinks=true,
    linkcolor=blue,
    citecolor=blue
}
\newtheorem{Th}{\scshape  Theorem}[section]
\newtheorem{Lem}[Th]{\scshape  Lemma}
\newtheorem{Cor}[Th]{\scshape Corollary}

\newtheorem{Prop}[Th]{\scshape Proposition}

\title{Line graphs  with the largest eigenvalue multiplicity}

\author{Wenhao Zhen\thanks{Corresponding author, E-mail address:zhenwenhao994@163.com.}   
 \ \  Dein Wong\thanks{Corresponding author, E-mail address:wongdein@163.com.    Supported by the National Natural Science Foundation of China
(No.12371025).}, \ \ Songnian Xu 
\\  {\small  \it  School of Mathematics, China University of Mining and Technology, Xuzhou, China} }

\date{}

\begin{document}
	\maketitle
	
\noindent \textbf{Abstract: }
{For a connected graph $G$, we denote by $L(G)$, $m_{G}(\lambda)$, $c(G)$ and $p(G)$ the line graph of $G$, the eigenvalue multiplicity of $\lambda$ in $G$, the cyclomatic number and the number of pendant vertices in $G$, respectively.
In 2023, Yang et al. \cite{WL LT} proved that $m_{L(T)}(\lambda)\leq p(T)-1$ for any tree $T$ with $p(T)\geq 3$, and characterized all trees $T$ with $m_{L(T)}(\lambda) = p(T)-1$.
In 2024, Chang et al. \cite{-1 LG} proved that, if $G$ is not a cycle, then $m_{L(G)}(\lambda)\leq 2c(G)+p(G)-1$, and they characterized all graphs $G$ with $m_{L(G)}(-1) = 2c(G)+p(G)-1$.
The authors of \cite{-1 LG}  particularly stated that it seems somewhat difficult to characterize  the extremal graphs $G$ with $m_{L(G)}(\lambda)= 2c(G)+p(G)-1$ for an arbitrary eigenvalue $\lambda$ of $L(G)$. In this paper, we give this problem a complete solution.
}


\section{Introduction}
In this paper, all graphs are   simple, i.e., they have no loops nor multiple edges. For a graph $G$ with
  vertex set   $V(G)$ and edge set $E(G)$, the adjacency matrix $A(G)$ of   $G$ with $|V(G)|=n$  is an $n\times n$ symmetric matrix whose $(i,j)$ entry   is $1$ if there is an edge between vertex $i$ and vertex $j$, and $0$ otherwise.
The eigenvalues of $A(G)$ are called the  eigenvalues of   $G$. Denote by $m_{G}(\lambda)$ the  multiplicity of $\lambda$ as an eigenvalue of $A(G)$.
If $\lambda=0$, $m_{G}(\lambda)$ is also called the \emph{nullity} of $G$.
The \emph{line graph} $L(G)$ of   $G$ is defined as the graph whose vertices are the edges of $G$, and two vertices in $L(G)$ are adjacent if and only if the corresponding edges in $G$ share exactly one common vertex.
For $u\not=v\in V(G)$, we write $u \sim v$ to denote that $u,v$ are adjacent vertices.
Let $N_G(u) = \{v| v \sim u, v \in V(G)\}$ denote  {the set of neighbors} of   $u$ in   $G$.
The \emph{degree} of  $u$ in   $G$ is defined as the number of neighbors of $u$ in $G$, which is denoted by $d_G(u)$.
A vertex $u\in G$ with $d_G(u)\geq 3$ (resp., $d_G(u)=1$) is called a \emph{major} (resp., \emph{pendant}) vertex of $G$.
The number of pendant vertices in   $G$ is denoted by $p(G)$.
If $G$ is connected, we call $|E(G)|-|V(G)|+1$ the \emph{cyclomatic} number of  $G$ and denote it by $c(G)$.

Some publications have concentrated on bounding the multiplicity for an arbitrary eigenvalue of a graph in terms of some of its structural parameters.
In \cite{T1,T2}, the upper bound of $m_T{(\lambda)} $ in terms of the order of $T$ were investigated for a real number $\lambda$ , where $T$ is a tree.
In \cite{WL}, Wang et al. proved that if $G$ is not a cycle, then $m_{G}(\mu)\leq 2c(G)+p(G)-1$.
Furthermore, the extremal graphs $G$ with nullity $ 2c(G)+p(G)-1$ were characterized in \cite{N2} and \cite{WL2}, independently.
More recently, Zhang et al. \cite{YS} gave the complete characterization on the extremal graphs $G$ with
$m_G(\lambda)=2c(G)+p(G)-1$.

The eigenvalue multiplicity of   line graphs  also attracted some  attention.
Gutman and Sciriha \cite{Gutman line graph of tree 0} showed that the nullity of the line graph of a tree is at most one.
Li et al. \cite{nullity deep 1} demonstrated that the nullity of the line graph of a unicyclic graph is at most two.
For the line graph $L(T)$ of a tree $T$, Yang et al. \cite{WL LT} proved that $m_{L(T)}(\lambda)\leq p(T)-1$, where $\lambda$ is an eigenvalue of $L(T)$, and they characterized the extremal trees $T$ with $m_{L(T)}(\lambda) = p(T)-1$.
In 2024, Chang et al. \cite{-1 LG} proved that:
\begin{Prop}\label{Prop1} For any connected graph $G$, let $L(G)$ be its line graph. 
If $G$ is not a cycle, then  $m_{L(G)}(\lambda)\leq 2c(G)+p(G)-1.$ \end{Prop} 
The authors of \cite{-1 LG} have characterized the extremal graphs $G$ with $m_{L(G)}(-1) = 2c(G)+p(G)-1$.

Now a natural problem arises:
\vskip 2mm
\noindent {\sf Problem}\ \ {\it For a connected graph $G$ that is not a cycle with $m_{L(G)}(\lambda)= 2c(G)+p(G)-1$, what is possible value of $\lambda$?  Furthermore,  how about the extremal graphs $G$ with $m_{L(G)}(\lambda)= 2c(G)+p(G)-1$? }
 \vskip 2mm
\noindent

The authors of  \cite{-1 LG} particularly stated that this problem seems somewhat difficult. Indeed, the problem has been solved for two special cases: the case when $G$ is a tree   \cite{WL LT},  or the case when $\lambda=-1$ \cite{-1 LG}.
 However, for the general case when $G$ is not a tree or when $\lambda\not=-1$, the problem still keeps unknown.
The motivation of this paper is to give the problem a complete solution. 
The main result of this paper is as follows. 

\begin{Th}\label{main result}
Let $G\neq C_n$ be a connected graph, and $L(G)$ be its line graph. 
Then $m_{L(G)}(\lambda)=2c(G)+p(G)-1$ if and only if $\lambda$ and $G$ are one of the following forms: \\
(i) $\lambda=2\cos(\frac{i\pi}{m+1})$ where $i$ and $m+1$ are two co-prime positive integers with $1\leq i\leq m$; $G$ is a path such that $d(u,v)\equiv m (\bmod\ m+1)$ for any two distinct pendant vertices $v$ and $u$ of $G$. \\
(ii) $\lambda=2\cos(\frac{2k\pi}{2q+1})$ where $2q+1$ and $2k$ are two co-prime positive integers with $1\leq k\leq q$; $G$ is a tree such that $d(u,v)\equiv 2q (\bmod\ 2q+1)$ for any two distinct pendant vertices $v$ and $u$ of $G$ and $p(G)\geq 3$.\\
(iii) $\lambda=2\cos(\frac{i\pi}{m+1})$ where $i$ and $m+1$ are two co-prime positive integers with $1\leq i\leq m$; $G$ is obtained from a tree $T$, where $m_{L(T)}(\lambda)=p(T)-1$, by joining $c(G)$ cycles of order a multiple of $m+1$ (resp., $2(m+1)$) when $i$ is even (resp., odd) to $c(G)$ distinct pendent vertices of $T$, where $p(T)\geq c(G)$ and $1\leq c(G)\leq 2$. \\ 
(iv) $\lambda=2\cos(\frac{i\pi}{m+1})$ where $i$ and $m+1$ are two co-prime positive integers with $1\leq i\leq m$; $G$ is obtained from two cycles $C_1$ and $C_2$ by adding an edge joining $C_1$ and $C_2$, where $C_1$ and $C_2$ have order a multiple of $m+1$ (resp., $2(m+1)$) when $i$ is even (resp., odd).  \\
(v) $\lambda=2\cos(\frac{2k\pi}{2q+1})$ where $2q+1$ and $2k$ are two co-prime positive integers with $1\leq k\leq q$; $G$ is obtained from a tree $T$, where $m_{L(T)}(\lambda)=p(T)-1$, by joining $c(G)$ cycles of order a multiple of $2q+1$ to $c(G)$ distinct pendent vertices of $T$, where $p(T)\geq c(G)\geq 3 $. 
\end{Th}

\section{Preliminaries}
\begin{figure}
	\centering
	\includegraphics[width=0.7\linewidth]{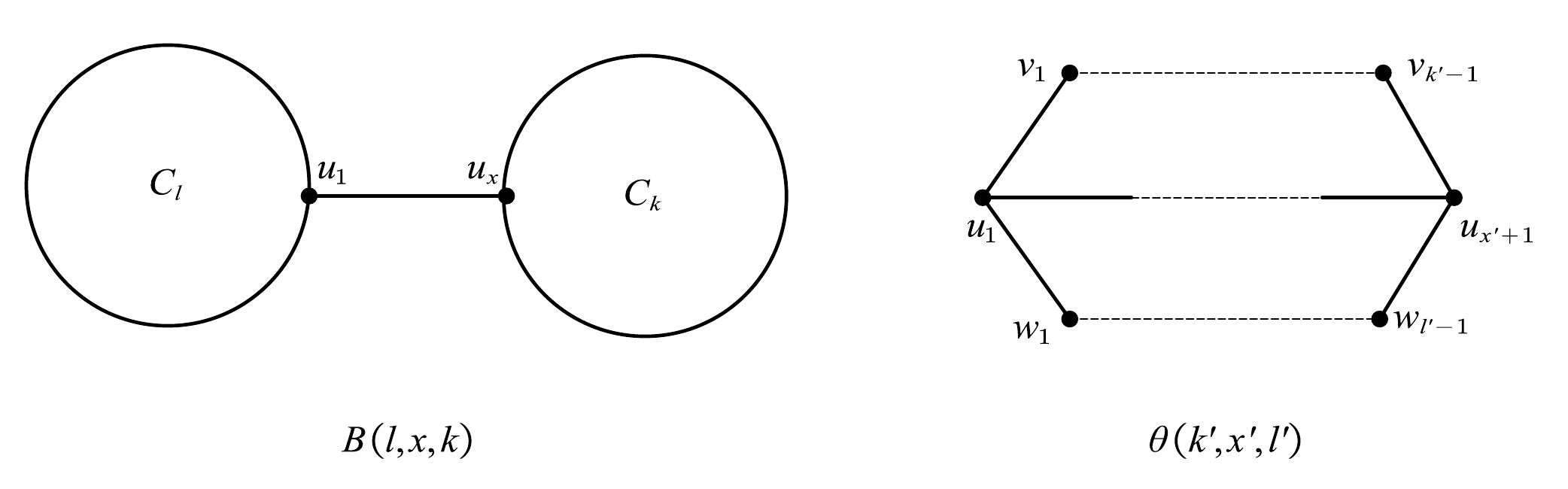}
	\caption{Graph $B(l,x,k)$ and $\theta (k',x',l',). $}
	\label{fig:1}
\end{figure}

For a   connected graph $G$ that is not a cycle, we will say that $L(G)$ is \emph{$\lambda$-optimal} if $m_{L(G)}(\lambda)=2c(G)+p(G)-1$.
The aim of this paper is to give a characterization for all $\lambda$-optimal line graphs. Firstly, we introduce some basic symbols and concepts.
In a connected graph $G$, a vertex $u$ is called a \emph{cutpoint} if the removal of $u$ and its incident edges results in a disconnected graph.
For a subset $W$ of $V(G)$, we denote by $G[W]$ and  $G-W$ the \emph{induced} subgraph of $G$ with vertex set $W$ and $V(G)\backslash W$, respectively.
Let $G$ be a graph, $x$ and $y$ be two vertices in $G$.
The \emph{distance} between $x$ and $y$ in $G$, denoted as $d_G(x,y)$, is defined as the length of a shortest path between them.
We denote by $P_n$ (resp., $C_n$)   a path (resp., cycle) with $n$ vertices.
Let $B(l,x,k)$ be the bicyclic graph (see the left side graph in Figure \ref{fig:1}) which is obtained from two vertex-disjoint cycles $C_l$ and $C_k$ by connecting one vertex of $C_l$ and one vertex of $C_k$ with a path $P_x$, where $l,k\geq 3,x\geq 1$ (if $x=1$, it is equivalent to identifying a vertex of $C_l$ with a vertex of $C_k$).
Let $\theta (k',x',l')$ be the bicyclic graph (see the right side graph in Figure \ref{fig:1}) which is a union of three internally disjoint paths $P_{k'+1}$, $P_{x'+1}$ and $P_{l'+1}$ respectively with common end vertices, where $k', x', l'\geq 1$ and at most one of them is 1.

The following lemma is well-known (see \cite{bro}):

\begin{Lem}\label{Lem PCE}
Let $m$ be a positive integer.\\
(i) The eigenvalues of $P_m$ are $2\cos \frac{i\pi }{m+1}$, $i=1,\ldots, m$.\\
(ii) The eigenvalues of $C_m$ are $2\cos \frac{2k\pi }{m}$, $k=0,\ldots, m-1$.
\end{Lem}

%
%

After   simple calculation  a lemma follows immediately from \ref{Lem PCE}.

\begin{Lem}\label{Lem lambda}
Suppose $\lambda =2\cos \frac{i\pi }{m+1}$, where $i$ and $m+1$  are two co-prime integers with $1\leq i\leq m$. \\
(i) $\lambda$ is an eigenvalue of $P_k$ if and only if $k\equiv m(\bmod\ {m+1})$.\\
(ii) If $i$ is even, then $\lambda$ is an eigenvalue of $C_k$ with $m_{C_k}(\lambda)=2$ if and only if $k$ is a multiple of $m+1$.\\
(iii) If $i$ is odd, then $\lambda$ is an eigenvalue of $C_k$ with $m_{C_k}(\lambda)=2$ if and only if $k$ is a multiple of $2(m+1)$. 
\end{Lem}

For a graph  $G$   with $n$ vertices, let $\mathbb{V}^{\lambda}_G = \{ x \in \mathbb{R}^n| A(G)x = \lambda x\}$ be the {characteristic subspace} of $A(G)$ with respect to  $\lambda$, where $\mathbb{R}^n$ is the $n$-dimensional column  {vector space} over the field of real numbers. 
The characteristic equation with respect to $\lambda$ at vertex $u\in V(G)$ is $\lambda x_u-\sum_{v\in N_G(u)}x_v=0$, where $x_u$ denotes the entry of $x$ at vector $x$. 
For a subset $U$  of $V(G)$, set $$\mathbb{Z}_G(U) = \{x\in \mathbb{R}^n| x_u = 0, \forall u \in U\}.$$ 
Then $\mathbb{Z}_G(U)$ is a subspace of $\mathbb{R}^n$ of dimension $n-|U|$.
If $\mathbb{V}^{\lambda}_G \cap \mathbb{Z}_G(U) = 0$, then $U$ is called a \emph{$\lambda$-annihilator} of $A(G)$.
The order of a   {$\lambda$-annihilator} of $A(G)$ has close relation with $m_G(\lambda)$.
\begin{Lem}\label{Lem LHZ} {\rm (\cite{Wang.LHZ}, Lemma 2.1)}
If $U$ is a $\lambda$-annihilator with respect to $A(G)$, then $m_G(\lambda)\leq|U|$.
\end{Lem}

In $\cite{Zhen}$, this result has been generalized.

\begin{Lem}\label{Lem TGLHZ} {\rm (\cite{Zhen}, Lemma 2.2)}
Let $Y\subseteq X$ be a subset  of $V(G)$ with $|X|-|Y|=m$.
Then $\dim(\mathbb{V}_G^{\lambda}\cap \mathbb{Z}_G(Y))\leq \dim(\mathbb{V}_G^{\lambda}\cap \mathbb{Z}_G(X))+m$.
Furthermore, $m_G(\lambda)\leq dim(\mathbb{V}_G^{\lambda}\cap \mathbb{Z}_G(X))+|X|$.
\end{Lem}

\begin{Lem}\label{Lem ob} Let $G$ be a graph of order $n$ and let $U$ be a subset of $V(G)$.
If $x\in \mathbb{V}^{\lambda}_G \cap \mathbb{Z}_G(U)$, then $x|_{G-U}\in \mathbb{V}^{\lambda}_{G-U}$.
\end{Lem}

\begin{proof}
After appropriate labeling, the adjacency matrix $A$ of $G$ is given by $$A=\begin{pmatrix}
  A(G[U])&B \\
 B^T &A(G-U)
\end{pmatrix}. $$
Note that $x=\begin{pmatrix}
  0 \\
x|_{G-U}
\end{pmatrix}$ and $x\in\mathbb{V}^{\lambda}_G $, it follows from  $$Ax=\begin{pmatrix}
  A(G[U])&B \\
 B^T &A(G-U)
\end{pmatrix}\begin{pmatrix}
  0 \\
x|_{G-U}
\end{pmatrix}=\lambda\begin{pmatrix}
  0 \\
x|_{G-U}
\end{pmatrix}$$
 that $A(G-U)x|_{G-U}=\lambda x|_{G-U}$,   leading to the required result.
\end{proof}

Let $P_m$ be a path with $m\geq 1$, and let $H$ be a graph disjoint with $P_m$.
Let $G$ be a graph obtained from $H$ and $P_m$ by adding an edge joining a vertex $u$ of $H$ and a pendant vertex of $P_m$.
If $d_H(u)\geq 2$, we say $P_m$ is a pendant path of $G$ and $G-P_m$ is obtained from $G$ by an operation of \emph{path-deletion} (with respect to $P_m$).
A cycle $C$ in $G$ is called a \emph{pendant cycle} of $G$ if there is only one major vertex $u$ in $C$ and $d_G(u)=3$.

\begin{Lem}\label{Lem PD} Let $H$ be a graph obtained from $G$ by a path-deletion operation, then $$m_{L(G)}(\lambda)\leq m_{L(H)}(\lambda)+1.$$
\end{Lem}

\begin{proof}
Suppose $P_m$ is a pendant path of $G$ with $m\geq 1$, $v_1\sim v_2\sim \cdots \sim v_m$ are all vertices of $P_m$ and $H=G-P_m$. 
If $m=1$,  the conclusion holds naturally (by Interlacing Theorem).
Suppose $m\geq 2$.
Let $e_i$ be the edge connecting $v_i$ and $v_{i+1}$ for $1\leq i\leq m-1$ and $e_m$ be the edge connecting $v_m$ and a vertex of $H$ (see Figure \ref{fig:2}).
Let $v_{e_i}$ be the vertex of $L(G)$ corresponding to $e_i$, and $P'$ be the path of $L(G)$ with vertex set $\{v_{e_i}|1\leq i\leq m\}$.
Then graph $L(H)=L(G-P_m)$ is the graph $L(G)-P'$.

\begin{figure}
	\centering
	\includegraphics[width=0.7\linewidth]{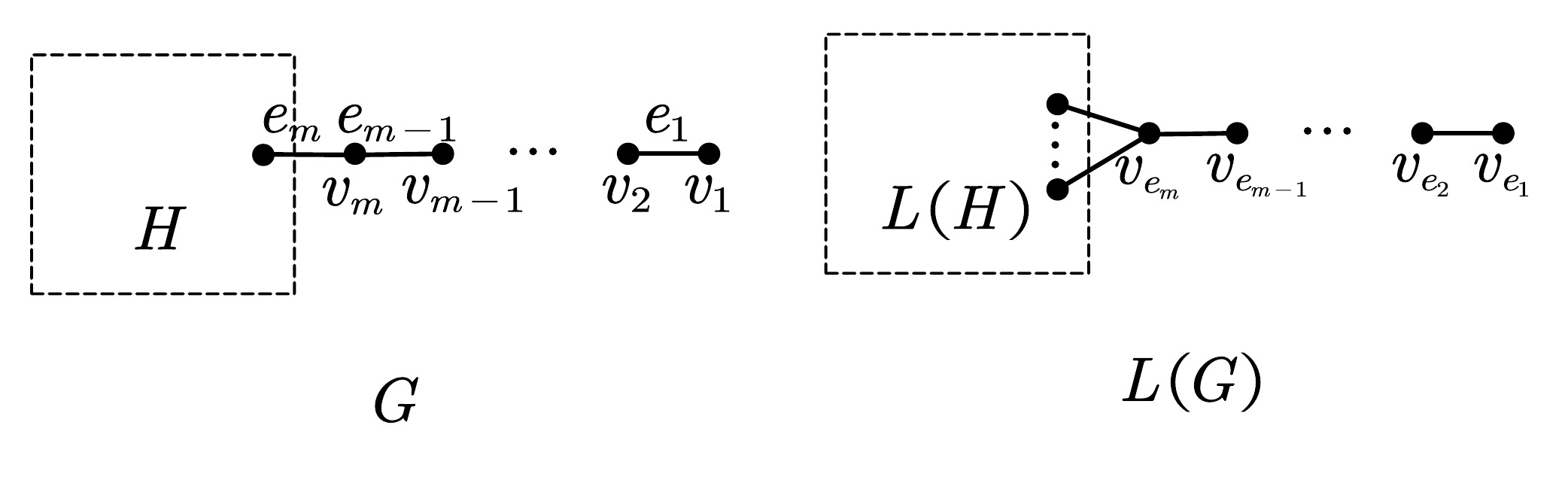}
	\caption{Graph $G$ and its line graph $L(G)$}
	\label{fig:2}
\end{figure}

Let $U=\{ v_{e_1} \}$.
If $x\in \mathbb{V}^{\lambda}_{L(G)} \cap \mathbb{Z}_{L(G)}(U)$, then $x_{v_{e_1}}=0$.
Since the unique neighbor of $v_{e_1}$ is $v_{e_2}$, by the characteristic equation at $v_{e_1}$, we have $\lambda x_{v_{e_1}}=x_{v_{e_2}}=0$.
Similarly, we have $x_{v_{e_1}}=x_{v_{e_2}}=\cdots =x_{v_{e_m}}=0$.
Thus, we have $x\in \mathbb{V}^{\lambda}_{L(G)} \cap \mathbb{Z}_{L(G)}(V(P'))$, which implies that $dim(\mathbb{V}^{\lambda}_{L(G)} \cap \mathbb{Z}_{L(G)}(U))\leq dim(\mathbb{V}^{\lambda}_{L(G)} \cap \mathbb{Z}_{L(G)}(V(P')))$.
By Lemma \ref{Lem ob}, we know $dim(\mathbb{V}^{\lambda}_{L(G)} \cap \mathbb{Z}_{L(G)}(V(P')))\leq m_{L(G)-P'}(\lambda)$.
By Lemma \ref{Lem TGLHZ}, we have
\begin{align*}
  m_{L(G)}(\lambda)&\leq dim(\mathbb{V}^{\lambda}_{L(G)} \cap \mathbb{Z}_{L(G)}(U))+1 \\
 &\leq dim(\mathbb{V}^{\lambda}_{L(G)} \cap \mathbb{Z}_{L(G)}(V(P')))+1\\
  &\leq m_{L(G)-P'}(\lambda)+1\\
   &=m_{L(H)}(\lambda)+1.
\end{align*}

\end{proof}

\begin{Lem}\label{Lem 2op}
Let $G$ be a connected graph that is not a cycle and $G'$ be obtained from $G$ by a path-deletion operation. If  $L(G)$ is $\lambda$-optimal and $G'$ is not a cycle,  then $L(G')$ is also $\lambda$-optimal.

\end{Lem}

\begin{proof}
Suppose $P$ is a pendent path of $G$ and $G'=G-P$, then $p(G)=p(G-P)+1$ and $c(G)=c(G-P)$.
By Lemma \ref{Lem PD}, we have
\begin{align*}
  m_{L(G-P)}(\lambda)&\geq m_{L(G)}(\lambda)-1\\
 &= 2c(G)+p(G)-2\\
  &= 2c(G-P)+p(G-P)-1.
\end{align*}
Note that $G'$ is not a cycle, we have $m_{L(G-P)}(\lambda)\leq 2c(G-P)+p(G-P)-1$ (by Proposition \ref{Prop1}). 
Hence, $m_{L(G-P)}(\lambda)= 2c(G-P)+p(G-P)-1$, which implies that $L(G')$ is $\lambda$-optimal.

\end{proof}

\begin{Lem} \label{john1}{\rm (\cite{john1}, Lemma 10)}\ Let $uv$ be a bridge of $G$. Suppose the induced subgraph on the two sides of the bridge are $G_u$ (containing $u$) and $G_v$ (containing $v$). If $\lambda$ is an eigenvalue of $G_u$ and $m(G_u, \lambda)=m(G_u-u, \lambda)+1$, then $m(G-v, \lambda)=m(G, \lambda)+1$.\end{Lem}
Applying Lemma 2.8, the following lemma is trivial.
\begin{Lem}\label{Lem Dm+1}
Suppose $\lambda =2\cos (\frac{i\pi}{m+1})$, where $i$ and $m+1$ are co-prime integers and $1\leq i\leq m$.
Let $G$ be a graph obtained from a graph $H$ and a disjoint path $P$ of order a multiple of $m+1$ by identifying one pendent vertex $v$ of $P$ with an arbitrary vertex of $H$.
Then $m_{G}(\lambda)= m_{G-P}(\lambda)$
\end{Lem}

\begin{proof}
Suppose    $u\in V(P)$ is adjacent to $v$. Then $uv$ is a bridge of $G$ and one side of the bridge containing $u$ is a path $P-v$ of order $m\   ({\rm mod}\  m+1$). Since $\lambda$ is a simple eigenvalue of $P-v$ and it is not an eigenvalue of $P-v-u$ (by Lemma \ref{Lem PCE}), we have $m(G-v, \lambda)=m(G, \lambda)+1$ (by Lemma \ref{john1}). As $G-v$ is the disjoint union of $P-v$ and $H-v$ and $m_{P-v}(\lambda)=1$, we have $m_{G}(\lambda)=m(G-v, \lambda)- 1=m_{G-P}(\lambda)$, as required.
  \end{proof}

\section{Characterization for all $\lambda$-optimal line graphs}

All $\lambda$-optimal line graphs of trees   have been characterized in \cite{WL LT}. Some notions related to line graphs $L(T)$ of trees should be introduced firstly.
A \emph{block} of a graph is a maximal induced subgraph without cutpoints.
The line graph $L(T)$ of a tree $T$ is a connected {block graph} in which each block is a clique and each cutpoint is in exactly two blocks.
In a connected graph $G$, a block with at least 3 vertices is called a \emph{major block}.
Let $B_{\mathfrak{L}}$ be the set of major blocks of $L(T)$.
A vertex of $L(T)$ shared by two different blocks is called an \emph{internal} vertex, otherwise it is an \emph{external} vertex.
Indeed, an external vertex of $L(T)$ corresponds to a pendant edge of $T$.
For an external vertex $v$ of $L(T)$, if there is $B\in B_{\mathfrak{L}}$ such that
$d(v,B)\leq d(v, B')$, $\forall B'\in B_{\mathfrak{L}}$, we call $v$ an  {external vertex of $B$} and call $B$ a  {major block of $v$}, where $d(v, B)=\min\{d(v, u)|u\in V(B)\}$.
If a major block $B$ of order $s$ has at least $s-1$ external vertices, $B$ is called an \emph{external} block.
If $L(T)$ has only one major block $B$, then $B$ is viewed as an external block of $L(T)$.
For $B_1, B_2\in B_{\mathfrak{L}}$, set
$d(B_1, B_2)=\min\{d(u_1, u_2)|u_i\in V(B_i)\}$, called the  {distance between $B_1$ and $B_2$}.

\begin{Th}\label{Th tree}{\rm (\cite{WL LT}, Theorem 1.2)}\ \
Let $T$ be a tree with $p\geq3$ pendant edges, and $L(T)$ be its line graph. Then $m_{L(T)}(\lambda)=p-1$ if and only if the following three conditions are all satisfied:\\
(i) There are two positive integers $q$ and $k$ with $1\leq k\leq q$ such that $\lambda=2\cos(\frac{2k\pi}{2q+1})$.\\
(ii) $d(v,B)+1\equiv q (\bmod\ 2q+1)$ for each external block $B$ and each external vertex $v$ of $L(T)$.\\
(iii) $d(B_1, B_2)\equiv 2q(\bmod\ 2q+1)$ for any two distinct major blocks $B_1$ and $B_2$ of $L(T)$.
\end{Th}

\noindent\textbf{{Remark.}}
  In condition (i), we can set $2q+1$ and $2k$ to be coprime positive integers. This is because if their greatest common divisor is not equal to $1$, after dividing them by this common divisor, they will still be an odd number and an even number respectively. 

For trees $T$ with $p(T)=2$, we can get the following result by Lemma \ref{Lem lambda} immediately.

\begin{Lem}\label{Th path}\ \
Let $T$ be a tree with $p(T)=2$, and $L(T)$ be its line graph. Then $m_{L(T)}(\lambda)=p(T)-1$ if and only if the following two conditions are all satisfied:\\
(i) There are two co-prime positive integers $i$ and $m+1$ with $1\leq i\leq m$ such that $\lambda=2\cos(\frac{i\pi}{m+1})$.\\
(ii) $d(u,v)\equiv m (\bmod\ m+1)$ for any two distinct pendant vertices $v$ and $u$ of $T$.
\end{Lem}

It is noteworthy that through a series of calculations, conditions (ii) and (iii) in Theorem \ref{Th tree} are found to be equivalent to $d(u,v)\equiv 2q-1 (\bmod\ 2q+1)$ for any two distinct external vertices $v$ and $u$ of $L(T)$, which leads to   the following corollary.

\begin{Cor}\label{Co tree}
Let $T$ be a tree, and $L(T)$ be its line graph. 
Then $m_{L(T)}(\lambda)=p-1$ if and only if $\lambda$ and $G$ are one of the following forms:\\
(i)$\lambda=2\cos(\frac{i\pi}{m+1})$ where $i$ and $m+1$ are two co-prime positive integers with $1\leq i\leq m$; $T$ is a path such that $d(u,v)\equiv m (\bmod\ m+1)$ for any two distinct pendant vertices $v$ and $u$ of $T$. \\
(ii)$\lambda=2\cos(\frac{2k\pi}{2q+1})$ where $2q+1$ and $2k$ are two co-prime positive integers with $1\leq k\leq q$; $T$ is a tree such that $d(u,v)\equiv 2q (\bmod\ 2q+1)$ for any two distinct pendant vertices $v$ and $u$ of $T$ and $p(T)\geq 3$.  
\end{Cor}

\begin{proof}
If $p(G)=2$, then by Lemma \ref{Th path}, we know this corollary holds. 
Suppose that $p(G)\geq 3$. 
Recall that an external vertex of $L(T)$ corresponds to a pendant edge of $T$. 
We only need to proof that conditions (ii) and (iii) in Theorem \ref{Th tree} hold if and only if $d(u,v)\equiv 2q-1 (\bmod\ 2q+1)$ for any two distinct external vertices $v$ and $u$ of $L(T)$.

Sufficiency: Let $T$ be a tree with $p\geq 3$ pendant vertices, $L(T)$ be its line graph and $d(u,v)\equiv 2q-1 (\bmod\ 2q+1)$ for any two distinct external vertices $v$ and $u$ of $L(T)$.
Suppose $B$ is a major block of $L(T)$, and $u_1$ is an external vertex of $L(T)$.
Then there must be a vertex $v_1\in V(B)$, such that $d(u_1,B)=d(u_1,v_1)$.
Let $v_2, v_3$ be two distinct vertices in $V(B-v_1)$.
Then there must be two external vertices $u_2,u_3$, such that $d(u_i,B)=d(u_i,v_i)$ for $i=2,3$ (see the left side graph in Figure \ref{fig:3}, and if $d(u_i,v_i)=0$ then $u_i=v_i$).
\begin{figure}
	\centering
	\includegraphics[width=0.9\linewidth]{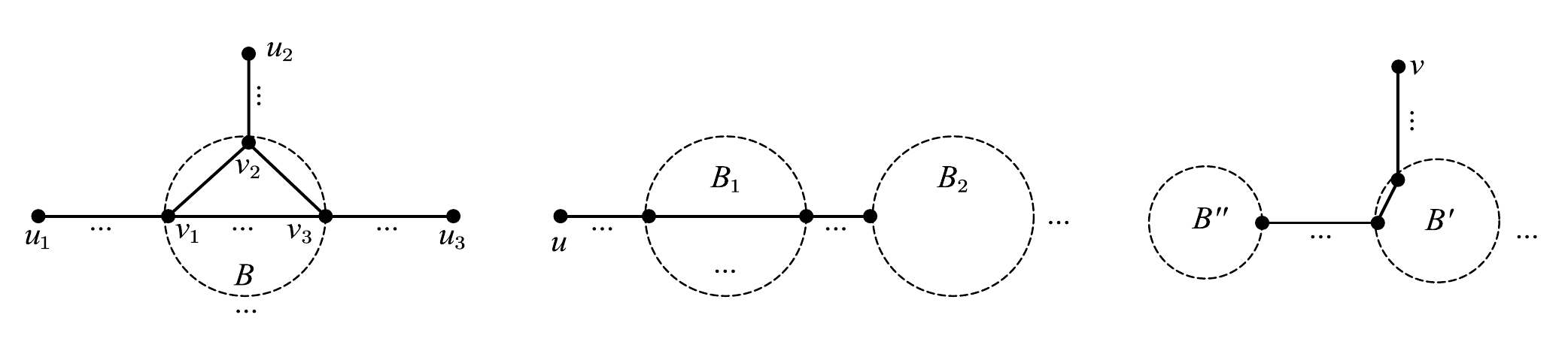}
	\caption{Graph $L(T)$.}
	\label{fig:3}
\end{figure}
Assume that $d(u_i,B)=a_i$ for $i=1,2,3$.
Then
\begin{align*}
d(u_1,u_2)&=a_1+a_2+1\equiv 2q-1(\bmod\ 2q+1)\\
d(u_1,u_3)&=a_1+a_3+1\equiv 2q-1(\bmod\ 2q+1)\\
d(u_2,u_3)&=a_2+a_3+1\equiv 2q-1(\bmod\ 2q+1).
\end{align*}
By calculating $d(u_1,u_2)+d(u_1,u_3)-d(u_2,u_3)$, we obtain $2a_1+1\equiv 2q-1(\bmod\ 2q+1)$), i.e., $a_1\equiv q-1(\bmod\ 2q+1)$.
Hence, we can see $d(v,B)+1\equiv q (\bmod\ 2q+1)$ for each major block $B$ and each external vertex $v$ of $L(T)$.
Note that any external block is a major block, one can see that condition (ii) in Theorem \ref{Th tree} holds.

Suppose $B_1$ and $B_2$ are two major blocks of $L(T)$.
Then there must be an external vertex $u$, such that $d(u,B_2)=d(u,B_1)+d(B_1,B_2)+1$ (see the middle graph in Figure \ref{fig:3}).
From the above discussion, we have $d(u,B_i)\equiv q-1 (\bmod\ 2q+1)$ for $i=1,2$.
Hence, $d(B_1,B_2)=d(u,B_2)-d(u,B_1)-1\equiv 2q(\bmod\ 2q+1)$, which implies that condition (iii) in Theorem \ref{Th tree} holds.

Necessity: Let $T$ be a tree with $p\geq 3$ pendant vertices and $L(T)$ satisfies conditions (ii) and (iii) in Theorem \ref{Th tree}. 
Suppose $v$ is an external vertex of $L(G)$, $B'$ is the major block of $v$. 
Firstly, we claim that $d(v,B')\equiv q-1(\bmod\ 2q-1)$. 
If $B'$ is an external block, then $d(v,B')\equiv q-1(\bmod\ 2q-1)$ by condition (ii). 
If $B'$ is not an external block, then there must be an external block $B''$ such that $d(v,B'')=d(B',B'')+d(v,B')+1$ (see the right side graph in Figure \ref{fig:3}). 
By conditions (ii) and (iii), we know $d(v,B'')=q-1(\bmod\ 2q+1)$ and $d(B',B'')=2q(\bmod\ 2q+1)$. 
Hence, we have $$d(v,B')=d(v,B'')-d(B',B'')-1\equiv q-1(\bmod\ 2q+1).$$ 

Let $v'$ and $v''$ be two distinct external vertices of $L(G)$, and let $B_3$ and $B_4$ be the major blocks of $v'$ and $v''$, respectively. 
If $B_3\neq B_4$, then $$d(v',v'')=d(v',B_3)+d(B_3,B_4)+d(v'',B_4)+2\equiv 2q-1(\bmod\ 2q+1). $$
If $B_3=B_4$, then $$d(v',v'')=d(v',B_3)+d(v'',B_3)+1\equiv 2q-1(\bmod\ 2q+1). $$

\end{proof}

The following lemma plays a crucial role in characterizing extremal graphs.

\begin{Lem}\label{Tool}

Let $G\neq C_n$ be a connected graph of order $n$ with $c(G)\geq1$ and $L(G)$ be its line graph.
Then $m_{L(G)}(\lambda)=2c(G)+p(G)-1$ if and only if $m_{L(G)}(\lambda)=m_{L(G-e)}(\lambda)+1$, $m_{L(G-e)}(\lambda)=2c(G-e)+p(G-e)-1$ and $p(G-e)=p(G)+1$,
where $e$ is an edge lies in a cycle of $G$ and adjacent to a major vertex.
\end{Lem}

\begin{proof}
Let $e$ be an edge lies in a cycle of $G$ and adjacent to a major vertex and $v_e$ be the vertex in $L(G)$ corresponding to $e$.
Clearly, $G-e$ is connected and $L(G-e)=L(G)-v_e$.
If $G-e$ is not a cycle, then by the proof of Theorem 3.2 in \cite{-1 LG}, we have $m_{L(G)}(\lambda)=2c(G)+p(G)-1$ if and only if $m_{L(G)}(\lambda)=m_{L(G-e)}(\lambda)+1$, $m_{L(G-e)}(\lambda)=2c(G-e)+p(G-e)-1$ and $p(G-e)=p(G)+1$.

If $G-e$ is a cycle, then $G$ is a bicyclic graph $\theta(k,2,l)$.
By the proof of Theorem 3.2 in \cite{-1 LG}, we have $m_{L(G-e)}(\lambda)=2$.
According to Lemma \ref{Lem PCE}, we have $|\lambda|< 2$.
Suppose $v_1\sim v_2\sim \cdots \sim v_n\sim v_1$ is the cycle $L(G-e)$ and $N_{L(G)}(v_e)=\{v_1,\ v_n,\ v_i,\ v_{i+1}\}$ (see the right side graph in Figure \ref{fig:4}).
\begin{figure}
	\centering
	\includegraphics[width=0.4\linewidth]{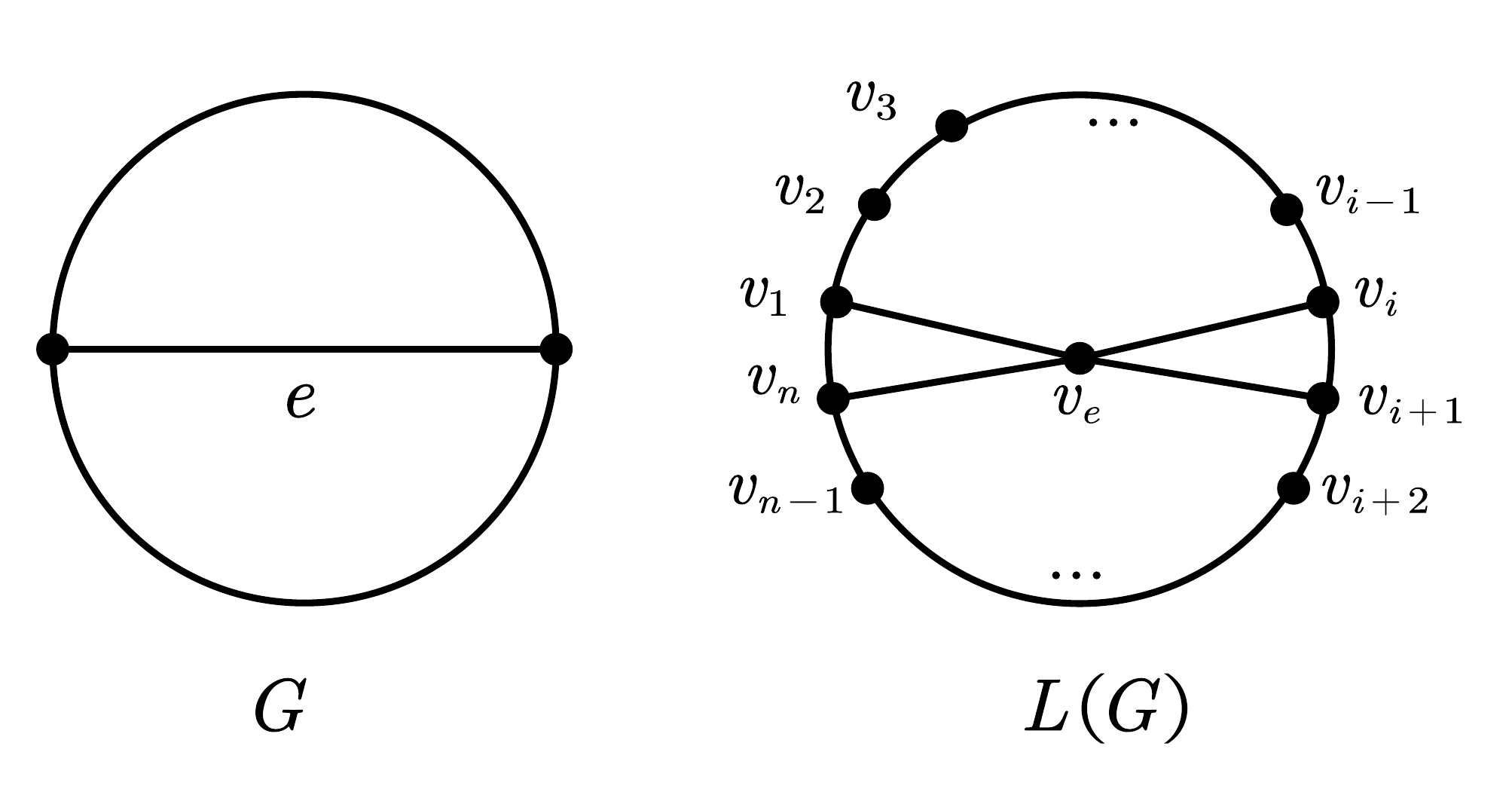}
	\caption{A bicyclic graph $G\cong \theta(k, 2, l)$ and its line graph $L(G)$.}
	\label{fig:4}
\end{figure}

Let $U=\{v_1, v_2\}$ and $x\in Z_{L(G)}(U)\cap V_{L(G)}^{\lambda}$.
Since $N_{L(G)}(v_2)=\{v_1,v_3\}$, it follows from $\lambda x_{v_2}=x_{v_1}+ x_{v_3}$ that $x_{v_3}=0$.
Similarly, we have $v_1=v_2=\cdots=v_i=0$.
Since $N_{L(G)}(v_i)=\{v_e,v_{i+1},v_{i-1}\}$, we have $x_{v_e}+x_{v_{i+1}}=0$.
Assume $x_{v_{i+1}}=a$ and $x_{v_{e}}=-a$.
Note that $N_{L(G)}(v_1)=\{v_2,v_n,v_{e}\}$, we have $x_{v_e}+x_{v_{n}}=0$.
This implies that $x_{v_{n}}=a$.
Moreover, since  $N_{L(G)}(v_e)=\{v_1,v_n,v_{i-1},v_i\}$, we have $-\lambda a=2a$.
If $a\neq 0$, we have $\lambda=-2$, which contradicts $|\lambda|<2$.
Hence, $a=0$ and we can get $x_{v}=0$ for any $v\in V(L(G))$ (by similar discussion as above).
By Lemma \ref{Lem LHZ}, we have $$m_{L(G)}(\lambda)\leq 2<2c(G)+p(G)-1,$$ a contradiction.
\end{proof}

\begin{Cor}\label{Co lambda}
Let $G\neq C_n$ be a connected graph and $L(G)$ be its line graph. If $L(G)$ is $\lambda$-optimal, then:\\
(i) There are two co-prime positive integers $i$ and $m+1$ with $1\leq i\leq m$ such that $\lambda=2\cos(\frac{i\pi}{m+1})$.\\
(ii) If $c(G)\geq 1$, then any two major vertices lies in a cycle are not adjacent in $G$.
\end{Cor}
Following from this corollary, we can deduce what is the possible value of an eigenvalue of a $\lambda$-optimal line graph.

\begin{proof}
(i) We proceed by induction on the cyclomatic number $c(G)$ of $G$.
If $c(G)=0$, then $G$ is tree and the  assertion holds by Corollary \ref{Co tree}.
Suppose that the assertion holds for graphs (different from a cycle) with cyclomatic number less than $c(G)$ and $c(G)\geq 1$.
Let $e$ be an edge lies in a cycle of $G$ and adjacent to a major vertex.
By Lemma \ref{Tool}, we have $m_{L(G-e)}(\lambda)=2c(G-e)+p(G-e)-1$.
Note that $G-e$ is still connected and $c(G-e)=c(G)-1$, by the induction hypothesis, we have $\lambda=2\cos(\frac{i\pi}{m+1})$ where $i$ and $m+1$ are co-prime positive integers.

(ii) Let $u$ be a major vertex in $G$ and $e$ be an edge lies in a cycle of $G$ and  incident to $u$.
By Lemma \ref{Tool}, we have $p(G-e)=p(G)+1$, which implies that the another endpoint of $e$ is not a major vertex.
\end{proof}

Now, we give a characterization on unicyclic graphs $G$ with $m_{L(G)}(\lambda)=2c(G)+p(G)-1=p(G)+1$.

\begin{Th}\label{Th U}
Let $G\neq C_n$ be a unicyclic graph and $\lambda=2\cos(\frac{i\pi}{m+1})$ where $i$ and $m+1$ are co-prime positive integers and $1\leq i\leq m$.
Then $L(G)$ is $\lambda$-optimal if and only if $G$ is obtained from a tree $T$, where $L(T)$ is $\lambda$-optimal, by joining a cycle of order a multiple of $m+1$ (resp., $2(m+1)$) when $i$ is even (resp., odd) to a pendent vertex of $T$.
\end{Th}

\begin{proof}
``If" part: Let $G\neq C_n$ be a unicyclic graph obtained from a tree $T$ that $L(T)$ is $\lambda$-optimal by joining a cycle of order a multiple of $m+1$ (resp., $2(m+1)$) when $i$ is even (resp., odd) to a pendent vertex of $T$. 
Let $e$ denote the added edge joining cycle $C$ and a pendent vertex of $T$, $v_e$ be the vertex of $L(G)$ corresponding  to $e$.
Clearly, $p(G)=p(T)-1$, $L(G)-v_e=L(G-e)=L(C)\cup L(T)$.
By Lemma \ref{Lem lambda}, we have $m_C(\lambda)=2$.
By the Interlacing Theorem, we have
\begin{align*}
 m_{L(G)}(\lambda)&\geq m_{L(G)-v_e}(\lambda)-1\\
 &= m_{L(C)}(\lambda)+m_{L(T)}(\lambda)-1\\
  &= 2+p(T)-1-1\\
  &=p(G)+1.
\end{align*}
Furthermore, $m_{L(G)}(\lambda)\leq 2c(G)+p(G)-1=p(G)+1$.
Hence, $m_{L(G)}(\lambda)=p(G)+1$.

``Only if" part: Let $G\neq C_n$ be a unicyclic graph that $L(G)$ is $\lambda$-optimal and let $C$ be the unique cycle in $G$.

$\mathbf{Claim 1.}$ $C$ has $\lambda$ as its eigenvalue.

After a series of path-deletion operations to $G$ (if necessary), we get a graph $G'$ that obtained from $C$ and a disjoint path $P$ by joining a pendant vertex of $P$ to a vertex of $C$.
Then $P$ is a pendant path of $G$.
By Lemma \ref{Lem 2op}, we have $m_{L(G')}(\lambda)=2c(G')+p(G')-1=2$.
By Lemma \ref{Lem PD} , we know
\begin{align*}
 m_{L(G')}(\lambda)&\leq m_{L(G-P)}(\lambda)+1 \\
 &= m_{L(C)}(\lambda)+1.
\end{align*}
Hence, $m_{C}(\lambda)=m_{L(C)}(\lambda)\geq 1$.
Moreover, by $|\lambda|<2$, we have $m_{C}(\lambda)=2$ and $C$ has order a multiple of $m+1$ (resp., $2(m+1)$) when $i$ is even (resp., odd) (By Lemma \ref{Lem lambda}). 

$\mathbf{Claim 2.}$ $C$ has only one major vertex.

Suppose $G$ has at least two major vertices. 
After a series of path-deletion operations to $G$ (if necessary), we get a graph $G''$ that is obtained from $C$ by joining a pendant vertex of two distinct paths respectively to two non-adjacent vertices of $C$ (see Figure \ref{fig:5}). 
The labeling of the vertices in $L(G'')$ is shown on the right side of Figure \ref{fig:5}. 
By Lemma \ref{Lem 2op}, we have $m_{L(G'')}(\lambda)=3$. 
\begin{figure}
	\centering
	\includegraphics[width=0.6\linewidth]{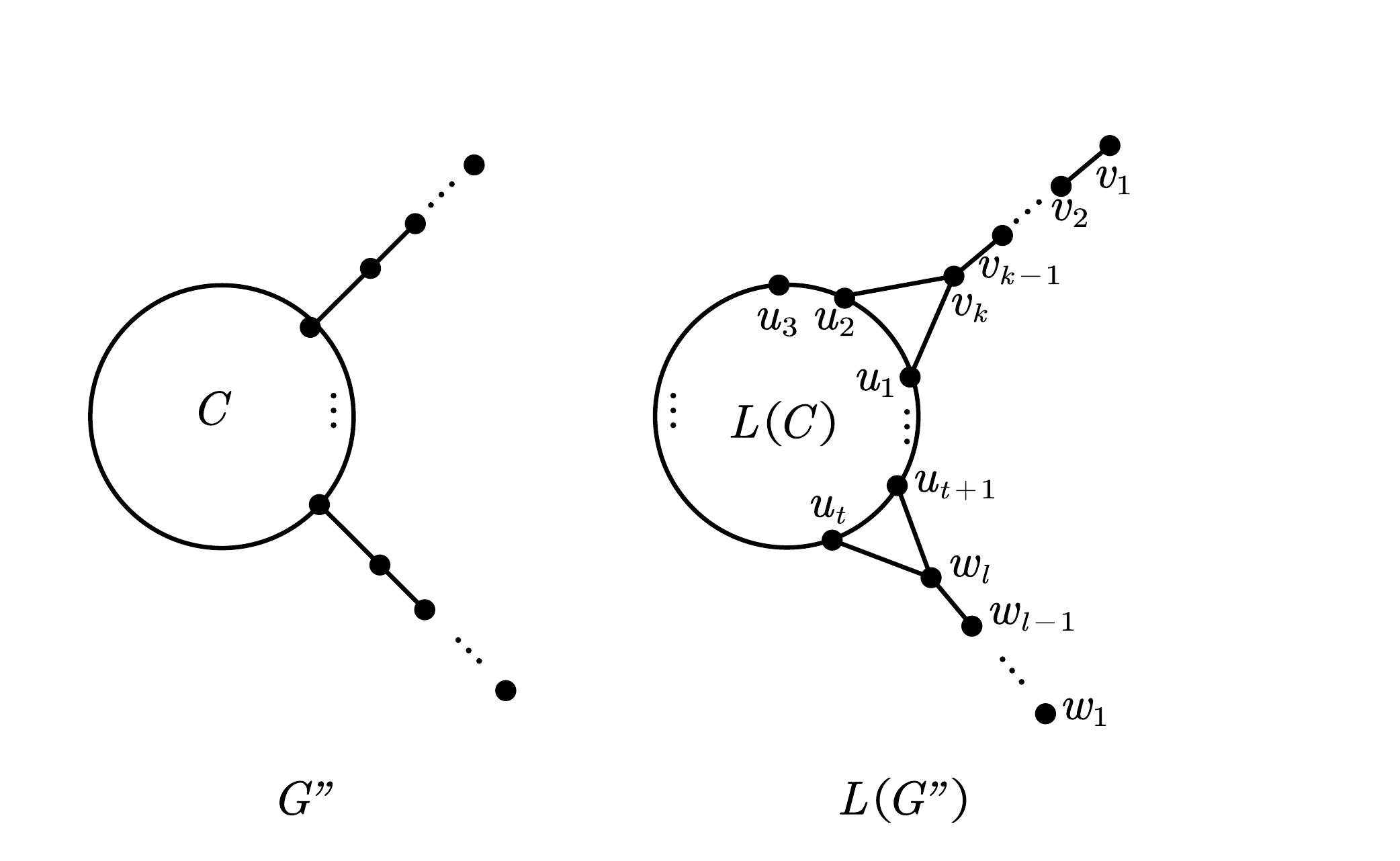}
	\caption{Graph $G''$ and its line graph $L(G'')$.}
	\label{fig:5}
\end{figure}

Let $U=\{v_1,u_1\}$ and $x\in \mathbb{Z}_{L(G)}(U)\cap \mathbb{V}_{L(G)}^{\lambda}$.
Since the only neighbor of $v_1$ is $v_2$, we have $x_{v_2}=\lambda x_{v_1}=0$.
Similarly, we know $x_{v_j}=0$ for $1\leq j\leq k$.
From the characteristic equation at vertex $v_k$, we have $$x_{u_2}=\lambda x_{v_k}-x_{v_{k-1}}-x_{u_1}=0.$$
From the characteristic equation at vertex $u_2$, we can see $x_{u_3}=0$.
Similarly, one can see $x|_{L(C)}=0$.
From the characteristic equation at vertex $u_t$, we have $x_{w_l}=0$.
Similarly, we have $x_{w_j}=0$ for $1\leq j\leq l$.
Finally, we have $x=0$.
By Lemma \ref{Lem LHZ}, we have $m_{L(G'')}(\lambda)\leq 2$, a contradiction.

Let $u$ be the only major vertex of $C$.

$\mathbf{Claim 3.}$ $d_G(u)=3$.

Suppose $d(u)\geq 4$.
Let $e_1$ be an edge on $C$ adjacent with $u$ and let $u_1$ be the another endpoint of $e_1$. 
After a series of path-deletion operations to $G$ (if necessary), we get a graph $G'''$ that is obtained from $C$ by joining a pendant vertex of two distinct paths respectively to $u$ (see Figure \ref{fig:5}). 
Let $u_2$ and $u_3$ be the pendant vertices in $G'''$. 
By Lemma \ref{Lem 2op}, we have ${L(G''')}$ is $\lambda$-optimal.  
\begin{figure}
	\centering
	\includegraphics[width=0.6\linewidth]{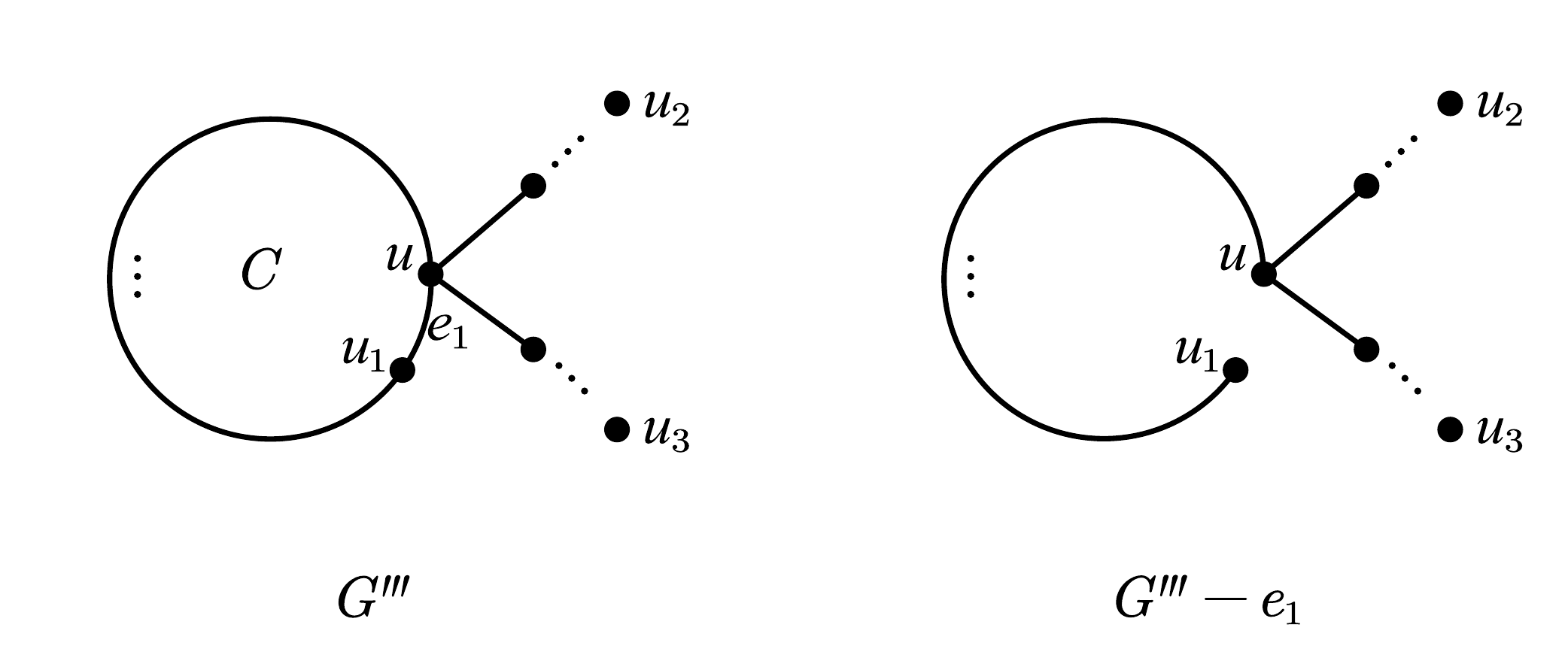}
	\caption{Graph $G'''$ and $G'''-e_1$.}
	\label{fig:6}
\end{figure}

By Lemma \ref{Tool}, we have $L(G'''-e_1)$ is $\lambda$-optimal. 
Note that $G'''-e_1$ is a tree, we have $d_{G'''-{e_1}}(u_i,u_j)\equiv m(\bmod\ m+1)$ for $i,j\in\{1,2,3\}$ and $i\neq j$ (By Corollary \ref{Co tree}).  
By Claim 1, we know that the order of $C$ is a multiple of $m+1$, which implies that $$d_{G'''-{e_1}}(u_1,u)\equiv m (\bmod\ m+1).$$ 
Hence, we have 
\begin{align*}
 d_{G'''-{e_1}}(u_2,u_3)&=d_{G'''-{e_1}}(u_2,u)+d_{G'''-{e_1}}(u_3,u)\\ 
 &=d_{G'''-{e_1}}(u_2,u_1)+d_{G'''-{e_1}}(u_3,u_1)-2d_{G'''-{e_1}}(u_1,u)\\
 &\equiv 0 (\bmod\ m+1).
\end{align*}
This contradicts with $d_{G'''-{e_1}}(u_2,u_3)\equiv m(\bmod\ m+1)$ and $m\geq 1$.

Now, we know $G$ is obtained from a tree $T$ by joining $C$ to a vertex of $T$, where $T=G-C$. 
Let $y$ be the only neighbor of $u$ in $T$. 

$\mathbf{Claim 4.}$ $d_T(y)=1$. 

The proof of Claim 4 is similar to that of Claim 3, so it is omitted here.

\begin{figure}
	\centering
	\includegraphics[width=0.7\linewidth]{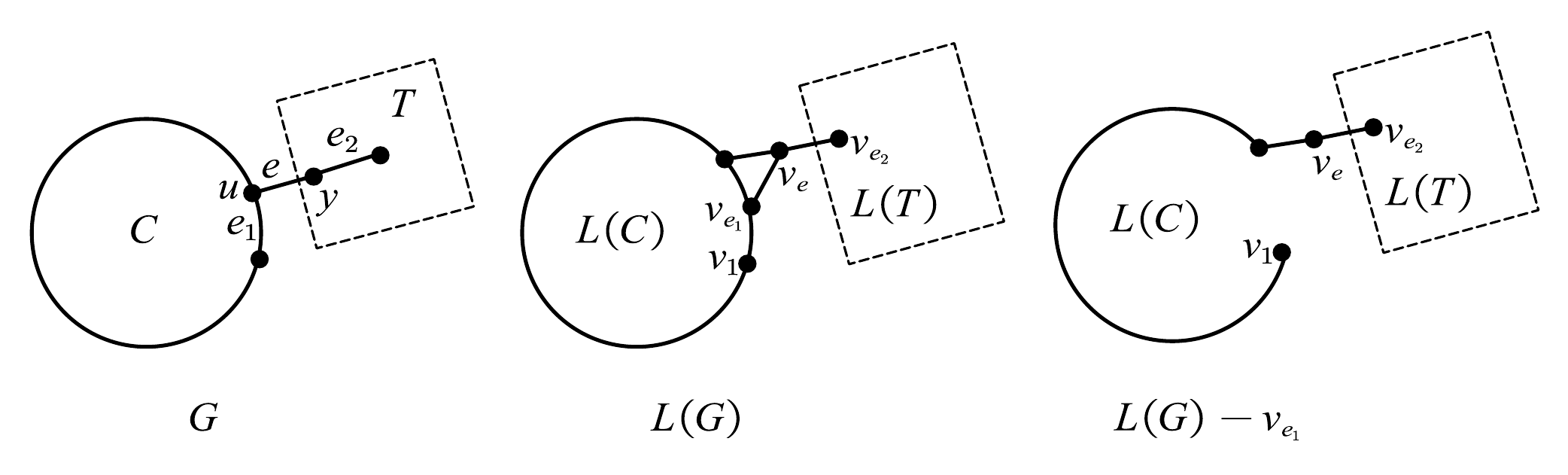}
	\caption{Graph $G$ and its line graph $L(G)$, and $L(G)-v_{e_1}$.}
	\label{fig:7}
\end{figure}

$\mathbf{Claim 5.}$ $L(T)$ is $\lambda$-optimal.

Let $e$ be the edge joining $y$ and $u$, and let $v_{e}$ be the vertex of $L(G)$ corresponding to $e$. 
Let $v_1$ be the 2-degree vertex on $L(C)$ that is adjacent to $v_{e_1}$ (see the graph in Figure \ref{fig:7}). 
Clearly, $$p(G-e_1)=p(T). $$
By Claim 1, we have $d_{L(G)-v_{e_1}}(v_1,v_e)\equiv m(\bmod\ m+1)$.
Recall that $$m_{L(G)-v_{e_1}}(\lambda)=p(G-e_1)-1.$$
By Lemma \ref{Lem Dm+1}, we have $$m_{L(T)}(\lambda)=m_{L(G)-v_{e_1}}(\lambda)=p(G-e_1)-1=p(T)-1.$$

According to above discussion, we know $G$ is obtained from a tree $T$, where $L(T)$ is $\lambda$-optimal, by joining a cycle of order a multiple of $m+1$ (resp., $2(m+1)$) when $i$ is even (resp., odd) to a pendent vertex of $T$.
\end{proof}

Next, we consider the extremal graphs of bicyclic graphs.

\begin{Th}\label{Th B}
Let $G$ be a bicyclic graph and $\lambda=2\cos(\frac{i\pi}{m+1})$ where $i$ and $m+1$ are co-prime positive integers and $1\leq i\leq m$.
Then $L(G)$ is $\lambda$-optimal if and only if $G$ is one of the following graphs:\\
(i) $G$ is obtained from a tree $T$, where $L(T)$ is $\lambda$-optimal, by joining two cycles of order a multiple of $m+1$ (resp., $2(m+1)$) when $i$ is even (resp., odd) to two distinct pendent vertices of $T$. \\
(ii) $G$ is obtained from two cycles $C_1$ and $C_2$ by adding an edge joining $C_1$ and $C_2$, where $C_1$ and $C_2$ have order a multiple of $m+1$ (resp., $2(m+1)$) when $i$ is even (resp., odd). 
\end{Th}

 \begin{proof}

``If" part: If $G$ is a bicyclic graph obtained from a tree $T$ that $L(T)$ is $\lambda$-optimal by joining two cycles $C_1$ and $C_2$ that have order a multiple of $m+1$ (resp., $2(m+1)$) when $i$ is even (resp., odd) to two distinct pendent vertex of $T$.
 Let $e_j$ denote the added edge joining cycle $C_j$ and a pendant vertex of $T$, and let $v_{e_j}$ be the vertex in $L(G)$ corresponding to $e_j$ for $j=1,2$, respectively.
Clearly, $p(G)=p(T)-2$, $L(G)-v_{e_1}-v_{e_2}=L(G-e_1-e_2)=L(C_1)\cup L(C_2)\cup L(T)$.
By Lemma \ref{Lem lambda}, we have $m_{C_j}(\lambda)=2$ for $j=1,2$.
By the Interlacing Theorem, we have
\begin{align*}
 m_{L(G)}(\lambda)&\geq m_{L(G)-v_{e_1}-v_{e_2}}(\lambda)-2\\
 &= m_{L(C_1)}(\lambda)+m_{L(C_2)}(\lambda)+m_{L(T)}(\lambda)-2\\
  &= 2+2+p(T)-1-2\\
  &=p(G)+3.
\end{align*}
Furthermore, $m_{L(G)}(\lambda)\leq 2c(G)+p(G)-1=p(G)+3$ (by Proposition \ref{Prop1}).
Hence, $m_{L(G)}(\lambda)=p(G)+3$. 

If $G$ is obtained from two cycles $C_1$ and $C_2$ by adding an edge joining $C_1$ and $C_2$, where $C_1$ and $C_2$ have order a multiple of $m+1$ (resp., $2(m+1)$) when $i$ is even (resp., odd).
Let $e$ be the edge joining $C_1$ and $C_2$, and let $v_e$ be the vertex of $L(G)$ corresponding to $e$.
Clearly, $L(G)-v_{e}=L(G-e)=L(C_1)\cup L(C_2)$.
By the Interlacing Theorem, we have
\begin{align*}
 m_{L(G)}(\lambda)&\geq m_{L(G)-v_{e}}(\lambda)-1\\
  &= m_{L(C_1)}(\lambda)+m_{L(C_2)}(\lambda)-1\\
  &= 2+2-1\\
  &=3.
\end{align*}
Since $m_{L(G)}(\lambda)\leq 2c(G)+p(G)-1=3$, we have $m_{L(G)}(\lambda)=3$.

``Only if" part:
Let $G$ be a bicyclic graph that $L(G)$ is $\lambda$-optimal.
Assume that two distinct cycles of $G$ are $C_1$ and $C_2$. Let $e_j$ be an edge lies in $C_j$ and adjacent to a major vertex for $j=1,2$ and let $v_{e_j}$ be the vertex of $L(G)$ corresponding to $e_j$ for $j=1,2$.

$\mathbf{Case 1. }$ $C_1$ and $C_2$ are vertex disjoint cycles.

By Lemma \ref{Tool}, we know $L(G-e_1)$ is $\lambda$-optimal.
Note that $G-e_1$ is a connected unicyclic graph, by Theorem \ref{Th U}, we know $C_2$ is a pendant cycle of $G-e_1$ and $C_2$ has order a multiple of $m+1$ (resp., $2(m+1)$) when $i$ is even (resp., odd).
Recall that $C_1$ and $C_2$ are vertex disjoint cycles, we have $C_2$ is also a pendant cycle of $G$.
Similarly, we can get that $C_1$ has order a multiple of $m+1$ (resp., $2(m+1)$) when $i$ is even (resp., odd) and $C_1$ is a pendant cycle of $G$.

If $G$ is not obtained from a tree $T$ by joining $C_1$ and $C_2$ to two distinct vertices of $T$.
Then $G$ is either obtained by adding an edge joining $C_1$ and $C_2$ or obtained by adding a vertex that adjacent to a vertex of $C_1$ and a vertex of $C_2$.

If $G$ is obtained by adding a vertex $u$ that adjacent to a vertex of $C_1$ and a vertex of $C_2$ (see Figure \ref{fig:8}).
\begin{figure}
	\centering
	\includegraphics[width=0.85\linewidth]{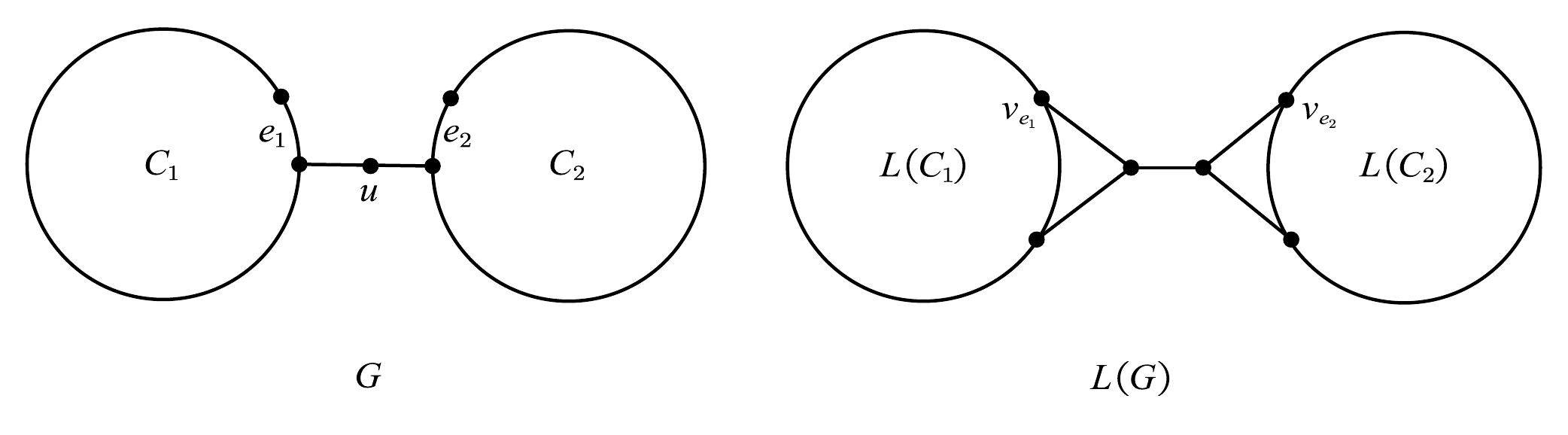}
	\caption{A bicyclic graph $G$ and its line graph $L(G)$.}
	\label{fig:8}
\end{figure}
Let $U=\{v_{e_1},v_{e_2}\}$ and $x\in \mathbb{Z}_{L(G)}(U)\cap \mathbb{V}_{L(G)}^{\lambda}$.
Recall that the order of $C_j$ is a multiple of $m+1$ for $j=1,2$.
Therefore, $L(G)-U$ is a path of order that is a multiple of $m+1$. 
By Lemma \ref{Lem ob}, we have $x|_{L(G)-U}\in \mathbb{V}^{\lambda}_{L(G)-U}$.
By Lemma \ref{Lem lambda}, we have $m_{L(G)-U}(\lambda)=0$, which implies that $x|_{L(G)-U}=0$.
Therefore, $x=0$.
By Lemma \ref{Lem LHZ}, we have $m_{L(G)}(\lambda)\leq 2$, a contradiction.

\begin{figure}
	\centering
	\includegraphics[width=0.85\linewidth]{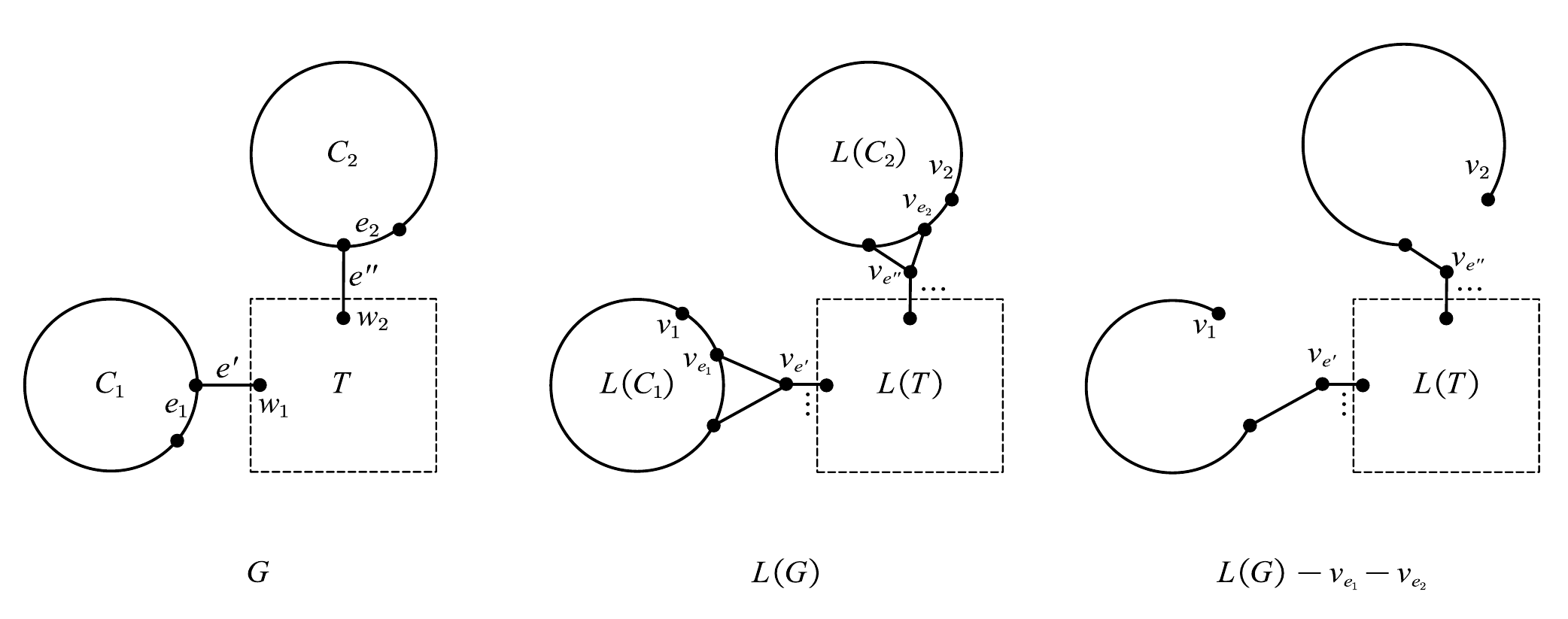}
	\caption{A bicyclic graph $G$ and its line graph $L(G)$, and $L(G)-v_{e_1}-v_{e_2}$.}
	\label{fig:9}
\end{figure}

If $G$ is obtained from a tree $T$ by joining $C_1$ and $C_2$ to two distinct vertices, say $w_j$ ($j = 1, 2$), of $T$  (see the left graph in Figure \ref{fig:9}).
By Theorem \ref{Th U}, we know $w_1$ is a pendant vertex of $G-e_2-C_1$.
Clearly, $w_1$ is also a pendant vertex of $T$.
Similarly, we can get $w_2$ is also a pendant vertex of $T$.
Then we only need to prove $m_{L(T)}(\lambda)=p(T)-1$.
The labeling of the vertices in $L(G)$ is shown
in Figure \ref{fig:9}.
Let $G'=L(G)-v_{e_1}-v_{e_2}=L(G-e_1-e_2)$.
Clearly, $p(G-e_1-e_2)=p(T).$
Recall that the order of $C_j$ is a multiple of $m+1$ for $j=1,2$, we have $d_{G'}(v_1,v_{e'})\equiv m(\bmod\ m+1)$ and $d_{G'}(v_2,v_{e''})\equiv m(\bmod\ m+1)$.
Applying Lemma \ref{Tool} two times, we have $$m_{G'}(\lambda)=p(G-e_1-e_2)-1.$$
By Lemma \ref{Lem Dm+1}, we have $$m_{L(T)}(\lambda)=m_{G'}(\lambda)=p(G-e_1-e_2)-1=p(T)-1.$$

$\mathbf{Case 2. }$ $C_1$ and $C_2$ share a common vertex.

If $C_1$ and $C_2$ have only one common vertex.
Recall that $C_2$ is a pendant cycle of $G-e_1$, we know that $C_2$ has only one major vertex.
Similarly we have $C_1$ has only one major vertex.
Hence $G\cong B(k,1,l)$.
Suppose $v_{e_1}\sim v_2\sim \cdots \sim v_k\sim v_{e_1}$ and $v_{e_2}\sim u_2\sim \cdots \sim u_l\sim v_{e_2}$ are $L(C_1)$ and $L(C_2)$, respectively (see Figure \ref{fig:10}).
Let $G''=L(G)-v_{e_1}-v_{e_2}$.
\begin{figure}
	\centering
	\includegraphics[width=0.85\linewidth]{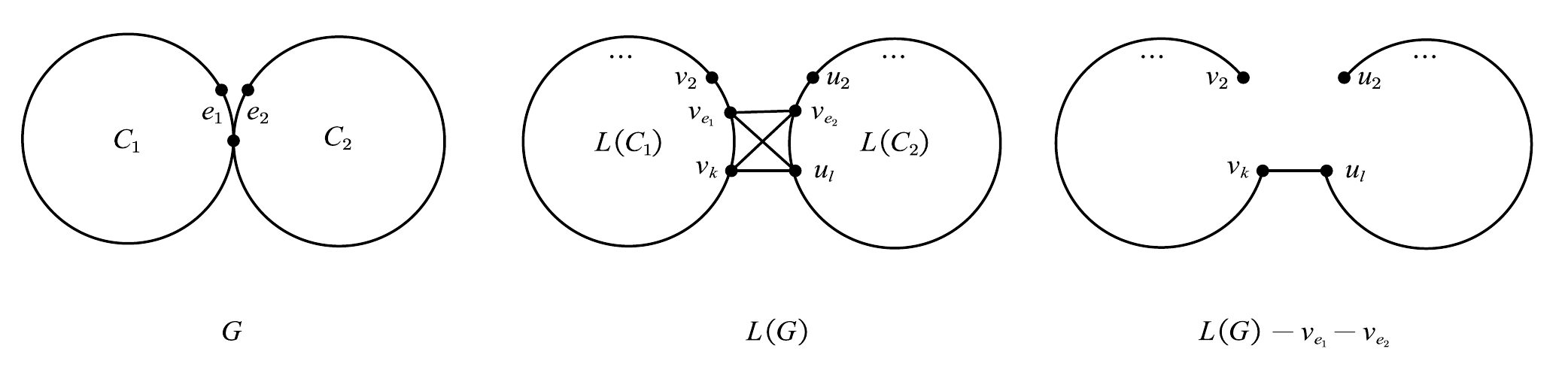}
	\caption{A bicyclic graph $G\cong B(k,1,l)$ and its line graph $L(G)$, and $L(G)-v_{e_1}-v_{e_2}$.}
	\label{fig:10}
\end{figure}

Let $U=\{v_{e_1},v_{e_2}\}$ and $x\in \mathbb{Z}_{L(G)}(U)\cap \mathbb{V}_{L(G)}^{\lambda}$.
By Lemma \ref{Lem ob}, we have $x|_{G''}\in \mathbb{V}^{\lambda}_{G''}$.
Similar to case 1, we know the order of $C_j$ is a multiple of $m+1$ for $j=1,2$. 
Therefore, we have $d_{G''}(v_2,u_2)\equiv m-2(\bmod\ m+1)$. 
Hence $G''$ is a path with order $d_{G''}(v_2,u_2)+1\equiv m-1(\bmod\ m+1)$. 
By Lemma \ref{Lem lambda}, we have $m_{G''}(\lambda)=0$, which implies that $x|_{L(G)-U}=0$.
Therefore, $x=0$.
By Lemma \ref{Lem LHZ}, we have $m_{L(G)}(\lambda)\leq 2$, a contradiction.

If $C_1$ and $C_2$ have at least two common vertices.
By Corollary \ref{Co lambda}, any two major vertices of $G$ lie in a cycle can not adjacent.
Recall that $C_1$ (resp., $C_2$) is a pendant cycle of $G-e_2$ (resp., $G-e_1$).
we have $G\cong \theta(k',x',l')$ which is a union
of three internally disjoint paths with common end vertices, where each path has length
at least 2 (see Figure \ref{fig:11}).
The labeling of the edges and vertices in $G$ is shown on Figure \ref{fig:11}.

\begin{figure}
	\centering
	\includegraphics[width=0.4\linewidth]{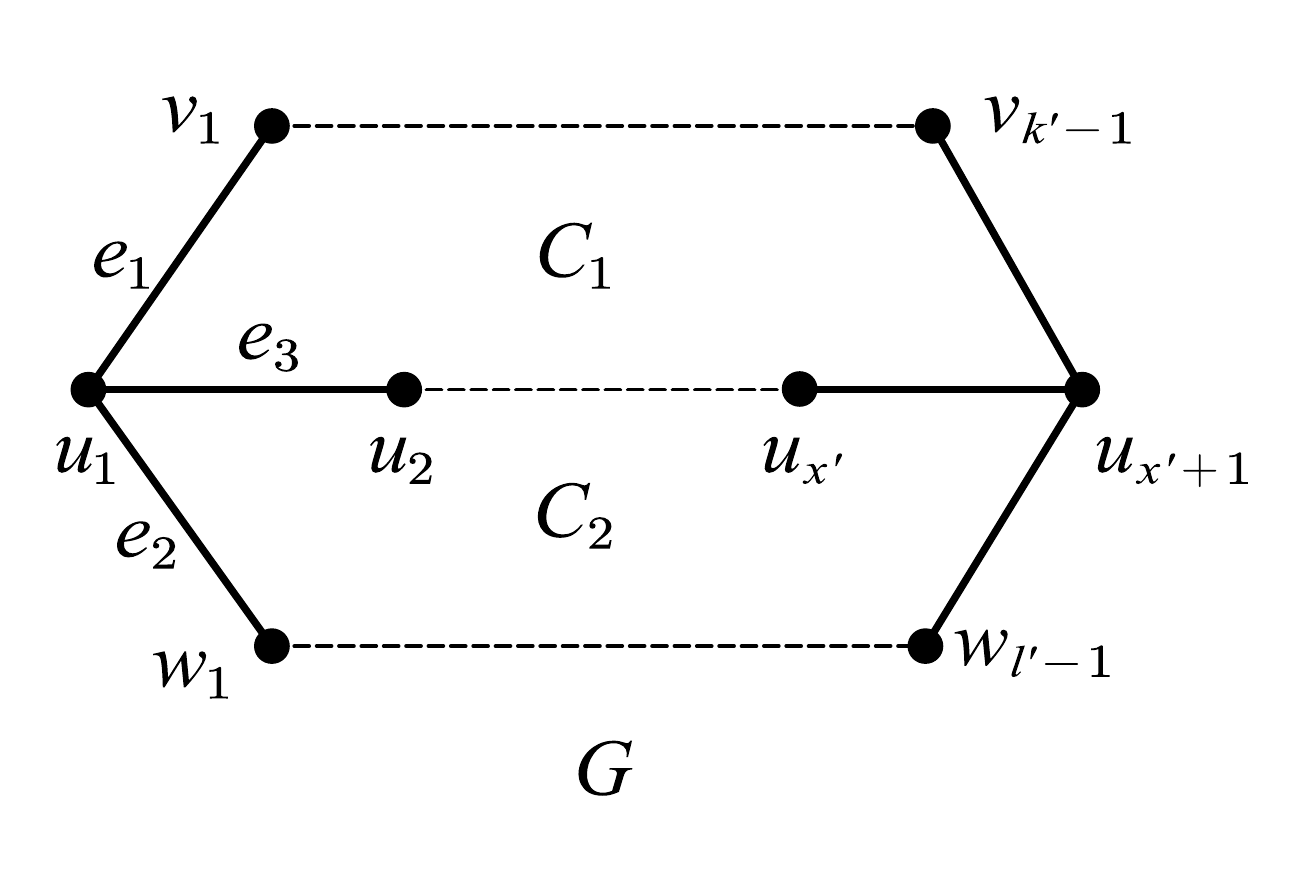}
	\caption{A bicyclic graph $G\cong \theta(k',x',l')$.}
	\label{fig:11}
\end{figure}

By Lemma \ref{Tool}, we have $L(G-e_1)$ is $\lambda$-optimal.
According to Theorem \ref{Th U} and Lemma \ref{Th path}, we have $d_{G-e_1}(v_1,v_{k'-1})\equiv m(\bmod\ m+1)$.
Similarly, we have $d_{G-e_2}(w_1,w_{l'-1})\equiv m(\bmod\ m+1)$.
Let $C$ be the cycle with vertex set $\{v_{1},v_{2},\dots ,v_{k'-1}, w_1, \cdots ,w_{l'-1}, u_1, u_{x'+1}\}$. 
Since $L(G-e_3)$ is $\lambda$-optimal, we know the order of $C$ is a multiple of $m+1$. 
Note that $|C|=d_{G-e_1}(v_1,v_{k'-1})+1+d_{G-e_2}(w_1,w_{l'-1})+1+2\equiv 2(\bmod\ m+1)$.
Hence, we have $m=1$, which implies that $\lambda=2\cos\frac{\pi}{2}=0$. 

By Lemma \ref{Lem lambda}, we know a cycle have $0$ as its eigenvalue if and only if the order of this cycle is a multiple of $4$.
We have $|C_j|\equiv0(\bmod\ 4)$ for $j=1,2$.
Since $L(G-e_3)$ is $\lambda$-optimal, we can get $x'$ is a odd number (by similar discussion as above). 
Hence, $|C|=|C_1|-x'+|C_2|-x'\equiv -2x'(\bmod\ 4)\equiv 2(\bmod\ 4)$, a contradiction.
 \end{proof}

Now, we are ready to give a proof for the main result.\\
\textbf{Proof of Theorem \ref{main result}
}
\begin{proof}
``If" part: 
If $\lambda$ and $G$ are one of the forms (i),(ii),(iii) or (iv), by Corollary \ref{Co tree}, Theorem \ref{Th U} and Theorem \ref{Th B}, we have $L(G)$ is $\lambda$-optimal. 

Let $\lambda=2\cos(\frac{2k\pi}{2q+1})$ where $2q+1$ and $2k$ are two co-prime positive integers with $1\leq k\leq q$. 
If $G$ is obtained from a tree $T$, where $m_{L(T)}=2c(T)+p(T)-1$, by joining $c(G)\geq 3$ cycles of order a multiple of $2q+1$ to $c(G)$ distinct pendent vertices of $T$, where $p(T)\geq c(G)\geq 3$. 
Let $e_j$ denote the added edge joining cycle $C_j$ and a pendent vertex of $T$, and let $v_{e_j}$ be the vertex of $L(G)$ corresponding to $e_j$ for $j = 1, 2, \cdots, c(G)$, respectively. 
Let $G'=G-e_1-\cdots -e_{c(G)}$. 
Clearly, $p(G)=p(T)-c(G)$, $L(G)-v_{e_1}-\cdots -v_{e_{c(G)}}=L(G')=C_1\cup \cdots \cup C_{c(G)}\cup L(T)$. 
By Lemma \ref{Lem lambda}, we have $m_{C_j}(\lambda)=2$ for $j=1,2,\dots ,c(G)$. 
By the Interlacing Theorem, we have 
\begin{align*}
 m_{L(G)}(\lambda)&\geq m_{G'}(\lambda)-c(G)\\
 &= m_{L(C_1)}(\lambda)+\cdots +m_{L(C_{c(G)})}(\lambda)+m_{L(T)}(\lambda)-c(G)\\
  &= 2c(G)+p(T)-1-c(G)\\
  &=2c(G)+p(G)-1.
\end{align*}
Furthermore, $m_{L(G)}(\lambda)\leq 2c(G)+p(G)-1$ (by Proposition \ref{Prop1}). 
Hence, $m_{L(G)}(\lambda)=2c(G)+p(G)-1$. 

``Only if" part: 
Let $G\neq C_n$ be a graph that $L(G)$ is $\lambda$-optimal. 
By Corollary \ref{Co lambda}, we know that there are two co-prime positive integers $i$ and $m+1$ with $1\leq i\leq m$ such that $\lambda=2\cos(\frac{i\pi}{m+1})$. 

Next, We proceed by induction on the cyclomatic number $c(G)$ of $G$ to prove $\lambda$ and $G$ must be one of the forms described in (i), (ii), (iii), (iv) and (v).  
If $c(G)=0,1,2$, then the assertion holds by Corollary \ref{Co tree}, Theorem \ref{Th U} and Theorem \ref{Th B}. 
Suppose that the assertion holds for graphs (different from a cycle) with cyclomatic number less than $c(G)$ and $c(G)\geq3$. 
First, we prove that $G$ is obtained from a tree $T$ by joining $c(G)$ cycles that have $\lambda$ as their eigenvalue to $c(G)$ distinct vertices. 
Choose any cycle $C$ of $G$, then there is a major vertex, say $u$, lies on $C$ since $G$ is not a cycle. 
Let $e$ be an edge lies on $C$ and adjacent to $u$, $x_e$ be the vertex of $L(G)$ corresponding to $e$. 
By Lemma \ref{Tool}, $L(G-e)$ is $\lambda$-optimal. 
Since $c(G-e)=c(G)-1\geq 2$, we know $G-e$ can't be a tree or unicyclic graph. 
If $G-e$ is obtained from two cycles $C_1$ and $C_2$ by adding an edge joining $C_1$ and $C_2$, then there must be another cycle $C'$ and an edge $e'$ on $C'$ that adjacent to a major vertex, such that $G-e'$ is a bicyclic graph that has two intersecting cycles. 
By Theorem \ref{Th B}, we know $L(G-e')$ is not $\lambda$-optimal, a contradiction. 
Then by the induction hypothesis, $G-e$ must be obtained from a tree $T'$, where $m_{L(T')}(\lambda)=p(T')-1$, by joining $c(G-e)$ cycles that have $\lambda$ as their eigenvalue to $c(G-e)$ distinct pendent vertices of $T'$. 
Note that $c(G)\geq 3$, there must be a cycle $C'$ disjoint with $C$ that have $\lambda$ as its eigenvalue in $G$. 
Clearly, $C'$ is a pendant cycle of $G$. 
If we choose $C'$ to begin the above discussion, then we have each cycle of $G$ is a pendant cycle of $G$ and $\lambda$ is a eigenvalue of each cycle.  
Then $G$ is obtained from a tree $T$ that $p(T)\geq c(G)$ by joining $c(G)$ cycles that have $\lambda$ as their eigenvalue to $c(G)$ distinct vertices, say $w_j$ ($j = 1, 2, . . . , c(G)$), of $T$. 
The labeling of the edges and vertices in $G$ and $L(G)$ is shown on Figure \ref{fig:12}. 

\begin{figure}
	\centering
	\includegraphics[width=0.95\linewidth]{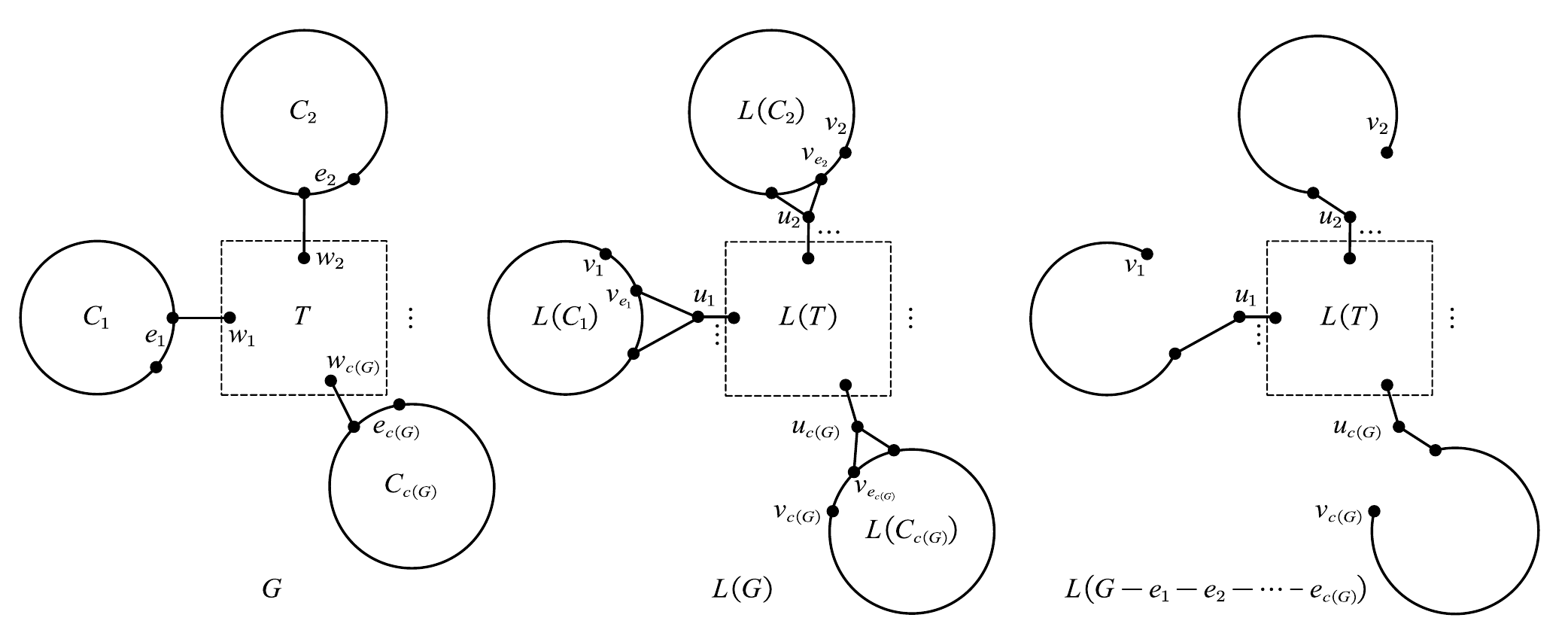}
	\caption{A graph $G$ and its line graph $L(G)$, and $L(G-e_1-\cdots -e_{c(G)})$.}
	\label{fig:12}
\end{figure}

By Lemma \ref{Tool}, we have $L(G-e_2-\cdots -e_{c(G)})$ is $\lambda$-optimal. 
By Theorem \ref{Th U}, we know $w_1$ is a pendant vertex of $G-e_2-\cdots -e_{c(G)}-C_1$. 
Clearly, $w_1$ is also a pendant vertex of $T$. 
Similarly, we have $w_j$ is a pendant vertex of $T$ for $j=1,2,\dots,c(G)$. 
Next, we prove $m_{L(T)}(\lambda)=p(T)-1$. 
Let $G'=L(G-e_1-e_2-\cdots e_{c(G)})$.  
Clearly, $$p(G-e_1-e_2-\cdots e_{c(G)})=p(T).$$ 
By Lemma \ref{Lem lambda}, we have $d_{G'}(v_j,u_j)\equiv m(\bmod\ m+1)$ for $j=1,2,\dots,c(G)$. 
Applying Lemma \ref{Tool} $c(G)$ times, we have $$m_{G'}(\lambda)=p(G-e_1-e_2-\cdots e_{c(G)})-1.$$  
By Lemma \ref{Lem Dm+1}, we have $$m_{L(T)}(\lambda)=m_{G'}(\lambda)=p(G-e_1-e_2-\cdots e_{c(G)})-1=p(T)-1.$$
Furthermore, since $p(T)\geq c(G)\geq 3$, we know there must be two positive integers $k$ and $q$ such that $i=2k$ and $m=2q$ (by Corollary \ref{Co tree}), which implies that $\lambda=2\cos(\frac{2k\pi}{2q+1})$. 
Moreover, by Lemma \ref{Lem lambda}, we know the order of each cycle of $G$ is a multiple of $2q+1$.

\end{proof}

\small {
	
}
\end{document}